\documentclass[12pt,a4paper]{article}
\usepackage[top=3cm, bottom=3cm, left=3cm, right=3cm]{geometry}
\usepackage[latin1]{inputenc}
\usepackage{csquotes}
\usepackage{amsmath}
\usepackage{amsthm}
\usepackage{amsfonts}
\usepackage{amssymb}
\usepackage{graphics}
\usepackage{float}
\usepackage[hang,flushmargin]{footmisc} 

\usepackage{amsmath,amsthm,amssymb,graphicx, multicol, array}
\usepackage{enumerate}
\usepackage{enumitem}
\usepackage{xcolor}
\usepackage{multicol}

\usepackage[pagebackref,bookmarks,colorlinks,breaklinks,linktoc=page]{hyperref}
\hypersetup{linkcolor=blue,citecolor=red,filecolor=blue,urlcolor=blue} 

\usepackage{tabto}

\numberwithin{equation}{section}

\usepackage{tikz}
\usepackage{pgfplots}

\usetikzlibrary{matrix}

\usepackage{mathrsfs}

\DeclareMathOperator{\li}{li}
\DeclareMathOperator{\ord}{ord}

\newtheorem{lem}{Lemma}[section]
\newtheorem{theorem}{Theorem}[section]
\newtheorem{lemma}{Lemma}[section]

\newcommand{\N}{\mathbb{N}}

\newcommand{\C}{\mathbb{C}}
\newcommand{\F}{\mathbb{F}}

\usepackage{fancyhdr}
\title{Primitive Root Conjecture in Arithmetic Progressions}
\date{}
\author{N. A. Carella}

\begin{document}
\maketitle


\begin{abstract}
Let \(x\geq 1\) be a sufficiently large number, and let $1 \leq a <q $ be a pair of integers such that $\gcd(a,q)=1$ and 
$q=O(\log^cx)$ with $c\geq 0$ constant. This note proves that the counting function for the number of primes $p \in \{p=qn+a: n \geq1 \}$ with 
a fixed squarefree primitive root $u\ne \pm 1$ has the asymptotic formula $\pi_u(x,q,a)=\delta(u,q,a)\li(x) +O(x/\log^b x),$ 
where $\delta(u,q,a)>0$ is the density, and $b=b(c)>1$ is a constant. \let\thefootnote\relax\footnote{ \today \date{} 
}
\end{abstract}
\tableofcontents
\section{Introduction}\label{S9900-I}\hypertarget{S9900-I}
Let $a\geq1$ and $q>a$ be a pair of integers such that $\gcd(a,q)=1$, and let $p\geq2$ be a prime. The multiplicative order modulo $p$ of 
an integer $u\ne \pm 1$ is denoted by $ \ord_p(u)=p-1$. The density of the subset of primes in the arithmetic progression 
$\{p=qn+a:\gcd(a,q)=1 \text{ and } \ord_p(u)=p-1\}$ is defined by a real number $\delta(u,q,a) \geq 0$, and the corresponding primes 
counting function is defined by
\begin{equation}
\pi _{u}(x,q,a)=\#\{p=qn+a\leq x :\gcd(a,q)=1 \text{ and } \ord_p(u)=p-1\},
\end{equation}
where $x \geq 1$ is a sufficiently large number. This note considers the followings result. 

\begin{theorem}   \label{thm9900.100}\hypertarget{thm9900.100}  Let $x \geq 1$ be a sufficiently large number. Let $1 \leq a <q$ be integers such that $\gcd(a,q)=1$ and 
$q=O(\log^c x)$, with $c\geq 0$ constant. Then, the arithmetic progression $\{p=qn+a:n \geq 1\}$ has infinitely many primes $p\geq 2$ with 
a fixed squarefree primitive root $u \ne \pm 1$. In addition, the corresponding primes counting function has the asymptotic formula
	\begin{equation} \label{e7l3}
	\pi _{u}(x,q,a)= \delta(u,q,a) \li(x)+O\left( \frac{x}{\log^b x} \right),
	\end{equation}
	where \(\delta(u,q,a) \geq 0\) is the density constant depending on the fixed integers $u,q,a $, and $b=b(c)>1$ is a constant.
\end{theorem} 

The preliminary foundational results and notation are discussed in \hyperlink{S9955D}{Section} \ref{S9955D} to \hyperlink{S9900-ET}{Section} \ref{S9900-ET}. The asymptotic formula for the main term is evaluated in \hyperlink{S9900-MT}{Section} \ref{S9900-MT} and an estimate of the error term is computed in \hyperlink{S9900-ET}{Section} \ref{S9900-ET}. Lastly, the proof of \hyperlink{thm9900.100}{Theorem} \ref{thm9900.100} is presented in \hyperlink{S9900-T1}{Section} \ref{S9900-T1}.

\section{Representations of the Characteristic Function}\label{S9955D}\hypertarget{S9955D}
The \textit{multiplicative order} of an element in a finite field $\mathbb{F}_p$ is defined by $\ord_p u=\min\{k:u^k\equiv 1 \bmod p\}$. An element $u\ne \pm1,v^2$ is called a primitive root if $\ord_p u=p-1.$ The characteristic function \(\Psi :G\longrightarrow \{ 0, 1 \}\) of primitive elements is one of the standard analytic tools employed to investigate the various properties of primitive roots in cyclic groups \(G\). Many equivalent representations of the characteristic function $\Psi $ of primitive elements
are possible, a few are investigated here. 
\subsection{Divisor Dependent Characteristic Function}		
The divisor dependent characteristic function was developed about a century ago, see {\color{red}\cite[Theorem 496]{LE1927}}
, {\color{red}\cite[p.\; 258]{LN1997}}, et alia. This characteristic function detects the order of an element by means of the divisors of the totient $p-1$. The precise description is stated below.

\begin{lemma} \label{lem9955.200D} \hypertarget{lem9955.200D} Let \(p\geq 2\) be a prime and let \(\chi\) be a multiplicative character of order $\ord  \chi =d$. If \(u\in
	\mathbb{F}_p\) is a nonzero element, then
	\begin{equation}
		\Psi (u)=\frac{\varphi(p-1)}{p-1}\sum _{d\mid p-1} \frac{\mu(d)}{\varphi(q)}\sum _{\ord \chi =d} \chi(u)
		=\left \{
		\begin{array}{ll}
			1 & \text{ if } \ord_p (u)=p-1,  \\
			0 & \text{ if } \ord_p (u)\neq p-1, \\
		\end{array} \right .\nonumber
	\end{equation}
	where $\mu:\N\longrightarrow \{-1,0,1\}$ is the Mobius function.
\end{lemma}	

\subsection{Divisorfree Characteristic Function}	
A new \textit{divisors-free} representation of the characteristic function of primitive element is developed here. It detects the order \(\text{ord}_p
(u) \geq 1\) of the element \(u\in \mathbb{F}_p\) by means of the solutions of the equation \(\tau ^n-u=0\) in \(\mathbb{F}_p\), where
\(u,\tau\) are constants, and $n$ is a variable such that \(1\leq n<p-1, \gcd (n,p-1)=1\). 

\begin{lemma} \label{lem9955.200A} \hypertarget{lem9955.200A} Let \(p\geq 2\) be a prime and let \(\tau\) be a primitive root mod \(p\) and  let \(\psi \neq 1\) be a nonprincipal additive character of order $\ord  \psi =p$. If \(u\in
	\mathbb{F}_p\) is a nonzero element, then
	\begin{equation}
		\Psi (u)=\sum _{\gcd (n,p-1)=1} \frac{1}{p}\sum _{0\leq s\leq p-1} \psi \left ((\tau ^n-u)s\right)
		=\left \{
		\begin{array}{ll}
			1 & \text{ if } \ord_p (u)=p-1,  \\
			0 & \text{ if } \ord_p (u)\neq p-1. \\
		\end{array} \right .\nonumber
	\end{equation}
\end{lemma}	
\begin{proof}[\textbf{Proof}] Set the additive character $\psi(s) =e^{i 2\pi  as/p}\in \C$. As the index $n\in \mathscr{R}=\{n<p:\gcd(n,p-1)=1\}$ ranges over the integers relatively prime to $\varphi(p-1)=p-1$, the element $\tau ^n\in \F_p ^{\times}$ ranges over the primitive roots
	modulo $p$. Accordingly, the equation $a=\tau ^n- u=0$ has a unique solution $n\geq1$ if and only if the fixed element $u\in \F_p$ is a primitive root. This implies that the inner sum in 	
	\begin{equation}\label{eq9977FF.300DF}
		\sum_{\gcd (n,p-1)=1} \frac{1}{p}\sum _{0\leq s< p} e^{i 2\pi \frac{(\tau ^n-u)s}{p}}=
		\left \{\begin{array}{ll}
			1 & \text{ if } \ord_{p} (u)=p-1,  \\
			0 & \text{ if } \ord_{p} (u)\ne p-1. \\
		\end{array} \right.
	\end{equation} 
	collapses to $\sum _{0\leq s< p} e^{i 2\pi as/p}=\sum _{0\leq s< p} 1=p $. Otherwise, if the element $u\in \F_p$ is not a primitive root, then the equation $a=\tau ^n- u=0$ has no solution $n\geq1$, and the inner sum in \eqref{eq9977FF.300DF} collapses to $\sum _{0\leq s< p} e^{i 2\pi as/p}=0$,
	this follows from the geometric series formula $\sum_{0\leq n\leq  N-1} w^n =(w^N-1)/(w-1)$, where $w=e^{i 2\pi a/p}\ne1$ and $N=p$. 
	This completes the verification.	 
\end{proof}

\section{Finite Summation Kernel and Gauss Sum}
An upper bound for some elementary exponential sums are provided in this section. 
\subsection{Finite Summation Kernel}
\begin{lemma}   \label{lem5555.400B}\hypertarget{lem5555.400B}  Let \(p\geq 2\) be large prime, and let $\omega=e^{i2 \pi/p} $ be a $p$th root of unity. Then,
\begin{equation}
		\sum_{1 \leq t\leq p-1}\Bigg | \sum_{\substack{1\leq n\leq p-1\\\gcd(n,p-1)=1}} \omega^{tn}  \Bigg  |\ll  p^{1+\delta} \log p ,
\end{equation}		where $\delta>0$ is a small number. 
\end{lemma} 

\begin{proof}[\textbf{Proof}] Use the inclusion exclusion principle to rewrite the exponential sum as
	\begin{eqnarray} \label{eq5555.400i}
		\sum_{\substack{1\leq n\leq p-1\\\gcd(n,p-1)=1}}\omega^{tn}&=& \sum_{n \leq p-1} \omega^{tn}  \sum_{\substack{d \mid p-1 \\ d \mid n}}\mu(d)  \nonumber \\
		&=& \sum_{d \mid p-1} \mu(d) \sum_{\substack{n \leq p-1 \\ d \mid n}} \omega^{tn}\nonumber \\
		& =&\sum_{d\mid p-1} \mu(d) \sum_{m \leq (p-1)/ d} \omega^{dtm} \\
		&=& \sum_{d \mid p-1} \mu(d) \frac{\omega^{dt}-\omega^{dt((p-1)/d+1)}}{1-\omega^{dt}} \nonumber.
	\end{eqnarray} 
	Now, the parameters are $p$ prime, $\omega=e^{i2 \pi/p}$, the integers $t \in [1, p-1]$ and $d \leq p-1<p$. This data implies that $\pi dt/p\ne k \pi $ with $k \in \mathbb{Z}$, so the sine function $\sin(\pi dt/p)\ne 0$ is well defined. Consequently, the absolute value satisfies
	\begin{equation}
		\left |\frac{\omega^{dt}-\omega^{dt((p-1)/d+1)}}{1-\omega^{dt}} \right |\leq 	\left | \frac{2}{\sin( \pi dt/ p)} \right |.
	\end{equation}
For each $d\mid p-1$, the map $t\longrightarrow z\equiv dt \bmod p$ is a permutation in the finite field $\F_p$. Thus, using standard manipulations, and $z/2 \leq \sin(z) <z$ for $0<|z|<\pi/2$, the last expression becomes
	\begin{eqnarray}
		\sum_{1 \leq t\leq p-1}\Bigg | 	\sum_{\substack{1\leq n\leq p-1\\\gcd(n,p-1)=1}}\omega^{tn}\Bigg  |&\leq&\sum_{d \mid p-1,} \sum_{1 \leq t\leq p-1}	\left | \frac{2}{\sin( \pi dt/ p)} \right |\\
		&\leq&\sum_{d \mid p-1,} \sum_{1 \leq z\leq p-1}	 \frac{2p}{ \pi z} \nonumber\\
		&\ll&  p^{1+\delta} \log p \nonumber,
	\end{eqnarray}
	where $\sum_{d \mid p-1}1=d(p-1)\ll p^{\delta}$ is the number of divisor in $p-1$ and $\delta>0$ is a small number. 
\end{proof}

\subsection{Gauss Sum}
Some elementary exponential sums estimates are provided in this section. 
\begin{lemma}   \label{lem1234A.150A}\hypertarget{lem1234A.150A}  
	{\normalfont (Gauss sums)} Let \(p\geq 2\) be a prime, let $\chi(t)=e^{i2 \pi t/p} $ and  $\psi(t)=e^{i2\pi  \tau^t/p}$ be a pair of characters. Then, the Gaussian sum has the upper bound
	\begin{equation} \label{eq3-355}
		\left |\sum_{1 \leq t \leq p-1}    \chi(t) \psi(t) \right | \leq 2 p^{1/2} \log p.\nonumber
	\end{equation}
	
\end{lemma}

\section{Estimates of Power Exponential Sums}
The estimate for the power sum with relatively prime index is based on the identity
\begin{equation}\label{eq9933Q.210c}
\frac{1}{p} \sum_{0 \leq t\leq p-1,}  \sum_{0 \leq s\leq p-1} \omega^{t(n-s)}f(s)=f(n),\end{equation}
where $\omega=e^{i2\pi/p}$ and $x  \leq p -1$.

\subsection{Power Exponential Sum with Relatively Prime Index}
\begin{theorem}  \label{thm9933Q.346}\hypertarget{thm9933Q.346}  Let \(p\geq 2\) be a large prime, and let $\tau $ be a primitive root modulo $p$. Then,
	\begin{equation}
		\sum_{\substack{1\leq n\leq p-1\\\gcd(n,p-1)=1}} e^{i2\pi b \tau^n/p} \ll  p^{1/2+\delta}(\log p)^2 \nonumber,
		\end{equation} 
		where $\delta>0$ is a small real number and the implied constant is independent of $b\ne0$. 	
\end{theorem}
\begin{proof}[\textbf{Proof}] Let $p$ be a large prime, and let $f(n)=e^{i 2 \pi b\tau^{n} /p}$, where $\tau$ is a primitive root modulo $p$. Start with the representation
	\begin{equation} \label{eq9933Q.346b}
		\sum_{\substack{1\leq n\leq p-1\\\gcd(n,p-1)=1}} e^{\frac{i2\pi b \tau^n}{p}}= \sum_{\substack{1\leq n\leq p-1\\\gcd(n,p-1)=1}}\frac{1}{p} \sum_{0 \leq t\leq p-1,}  \sum_{1 \leq s\leq p-1} \omega^{t(n-s)}e^{\frac{i2\pi b \tau^s}{p}} ,
	\end{equation}
see \eqref{eq9933Q.210c}. Use the inclusion exclusion principle to rewrite the exponential sum as
	\begin{equation}\label{eq9933Q.346d}
		\sum_{\substack{1\leq n\leq p-1\\\gcd(n,p-1)=1}} e^{ \frac{i2\pi b \tau^n}{p}} 
		= \sum_{1\leq  n \leq p-1}\frac{1}{p} \sum_{0 \leq t\leq p-1,}  \sum_{1 \leq s\leq p-1} \omega^{t(n-s)}e^{\frac{i2\pi b \tau^s}{p}} \sum_{\substack{d \mid p-1 \\ d \mid n}}\mu(d)   .
	\end{equation} 
	Now, observe that the term $t=0$ contributes $-\varphi(p-1)/p$, and rearranging it yield
	\begin{eqnarray}\label{eq9933Q.346f}
		&& \sum_{\substack{1\leq n\leq p-1\\\gcd(n,p-1)=1}} e^{ \frac{i2\pi b \tau^n}{p}} \\
		&=& \sum_{ n \leq p-1}\frac{1}{p} \sum_{1 \leq t\leq p-1,}  \sum_{1 \leq s\leq p-1} \omega^{t(n-s)}e^{\frac{i2\pi b \tau^s}{p}} \sum_{\substack{d \mid p-1 \\ d \mid n}}\mu(d) -\frac{\varphi(p-1)}{p} \nonumber \\
		&=&\frac{1}{p} \sum_{1 \leq t\leq p-1} \left ( \sum_{1 \leq s\leq p-1} \omega^{-ts}e^{\frac{i2\pi b \tau^s}{p}}\right )\left (\sum_{d \mid p-1} \mu(d) \sum_{\substack{n \leq p-1, \\ d \mid n}}   \omega^{tn} \right ) -\frac{\varphi(p-1)}{p} \nonumber.
	\end{eqnarray} 
	Taking absolute value, and applying \hyperlink{lem5555.400B}{Lemma} \ref{lem5555.400B}, and \hyperlink{lem1234A.150A}{Lemma} \ref{lem1234A.150A}, yield
	\begin{eqnarray} \label{eq9933Q.346h}
		&& \left | \sum_{\substack{1\leq n\leq p-1\\\gcd(n,p-1)=1}} e^{\frac{i2\pi b \tau^n}{p}} \right | \\
		&\leq&\frac{1}{p}  \sum_{1 \leq t\leq p-1} \left | \sum_{1 \leq s\leq p-1} \omega^{-ts}e^{i2\pi b \tau^{s}/p} \right | \cdot  \left |\sum_{d \mid p-1} \mu(d) \sum_{\substack{n \leq p-1, \\ d \mid n}}   \omega^{tn} \right | +\frac{\varphi(p-1)}{p}\nonumber \\
		&\ll&\frac{1}{p}  \sum_{1 \leq t\leq p-1} \left ( 2p^{1/2} \log p \right ) \cdot  \left |\sum_{d \mid p-1} \mu(d) \sum_{\substack{n \leq p-1, \\ d \mid n}}   \omega^{tn} \right |+\frac{\varphi(p-1)}{p}\nonumber\\
		&\ll&\frac{1}{p} \left ( 2p^{1/2} \log p \right ) \cdot   \sum_{1 \leq t\leq p-1} \left |\sum_{d \mid p-1} \mu(d) \sum_{\substack{n \leq p-1, \\ d \mid n}}   \omega^{tn} \right |+\frac{\varphi(p-1)}{p}\nonumber\\[.2cm]
		&\ll&\frac{1}{p} \left ( 2p^{1/2} \log p \right ) \cdot  \left ( 2p^{1+\delta} \log p \right ) \nonumber\\[.3cm]
		&\ll& p^{1/2+\delta} (\log p)^2 \nonumber,
	\end{eqnarray}
	where $\delta>0$ is a small number.
\end{proof}

A different approach to this result appears in {\color{red}\cite[Theorem 6]{FS2000}}, and related results are given in \cite{FS2001}, \cite{GM2005}, \cite{CC2009}, and {\color{red}\cite[Theorem 1]{GK2005}}. The upper bound given in \hyperlink{thm9933Q.346}{Theorem} \ref{thm9933Q.346} seems to be optimum. A different proof, which has a weaker upper bound, appears in {\color{red}\cite[Theorem 6]{FS2000}}, and related results are given in \cite{CC2009}, \cite{FS2001}, \cite{GK2005}, and {\color{red}\cite[Theorem 1]{GK2005}}.

\subsection{FFT of Power Exponential Sum with Relatively Prime Index} 
For any fixed $ 0 \ne b \in \mathbb{F}_p$, the map $ \tau^n \longrightarrow b \tau^n$ is one-to-one (permutation) in $\mathbb{F}_p$. Consequently, the subsets 
\begin{equation} \label{eq9933RPI.500b}
	\{ \tau^n: \gcd(n,p-1)=1 \}\quad \text { and } \quad  \{ b\tau^n: \gcd(n,p-1)=1 \} \subset \mathbb{F}_p
\end{equation} have the same cardinalities. As a direct consequence the exponential sums 
\begin{equation} \label{eq9933RPI.500d}
	\sum_{\substack{1\leq n\leq p-1\\\gcd(n,p-1)=1}}e^{i2\pi b \tau^n/p} \quad \text{ and } \quad \sum_{\substack{1\leq n\leq p-1\\\gcd(n,p-1)=1}} e^{i2\pi \tau^n/p},
\end{equation}
have the same upper bound up to an error term. An asymptotic relation for the finite Fourier transform (FFT) of the exponential sums (\ref{eq9933RPI.500d}) is provided here. 

\begin{theorem}   \label{thm9933ERP.220V}\hypertarget{thm9933ERP.220V}  Let \(p\geq 2\) be a large prime. If $\tau $ be a primitive root modulo $p$ and $a<x=o(p)$ is not a primitive root, then
	\begin{equation} 
	\widehat{V(a)}=	\sum_{1\leq b\leq  p-1}	 e^{-i2\pi \frac{ab}{p}}	\sum_{\substack{1\leq n\leq p-1\\\gcd(n,p-1)=1}} e^{\frac{i2\pi b \tau^n}{p}} =-  \sum_{\substack{1\leq n\leq p-1\\\gcd(n,p-1)=1}} e^{\frac{i2\pi  \tau^n}{p}} + O(p^{1/2+\delta} (\log p)^2)\nonumber,
	\end{equation} 
	where $\delta>0$ is a small number and the implied constant is independent of $ b \in [1, p-1]$. 	
\end{theorem} 
\begin{proof}[\textbf{Proof}] For $a\in[1,x]$ and $b\in[1,p-1]$, the exponential sum has the representation 
	\begin{eqnarray} \label{eq9933RPI.500f}
		V(b)&=& \sum_{\substack{1\leq n\leq p-1\\\gcd(n,p-1)=1}} e^{\frac{i2\pi b \tau^n}{p}} \\
		&=&\frac{1}{p} \sum_{1 \leq t\leq p-1} \left ( \sum_{1 \leq s\leq p-1} \omega^{-ts}e^{\frac{i2\pi b \tau^s}{p}}\right )\left (\sum_{d \mid p-1} \mu(d) \sum_{\substack{n \leq p-1, \\ d \mid n}}   \omega^{tn} \right ) -\frac{\varphi(p-1)}{p}\nonumber,
	\end{eqnarray} 
	confer equations \eqref{eq9933Q.346b}, \eqref{eq9933Q.346d} and \eqref{eq9933Q.346f} for more details. In particular, for $b=1$, 
	\begin{eqnarray} \label{eq9933RPI.500h}
		V(1)&=& 	\sum_{\substack{1\leq n\leq p-1\\\gcd(n,p-1)=1}}e^{\frac{i2\pi  \tau^n}{p}} \\
		&=& \frac{1}{p} \sum_{1 \leq t\leq p-1} \left ( \sum_{1 \leq s\leq p-1} \omega^{-ts}e^{\frac{i2\pi a \tau^s}{p}}\right )\left (\sum_{d \mid p-1} \mu(d) \sum_{\substack{n \leq p-1, \\ d \mid n}}   \omega^{tn} \right ) -\frac{\varphi(p-1)}{p}\nonumber,
	\end{eqnarray}
	respectively. Differencing (\ref{eq9933RPI.500f}) and (\ref{eq9933RPI.500h}) produces 
	\begin{eqnarray} \label{eq9933RPI.500i}
		V(b)-V(1)&= &	\sum_{\substack{1\leq n\leq p-1\\\gcd(n,p-1)=1}} e^{\frac{i2\pi b \tau^n}{p}} -\sum_{\substack{1\leq n\leq p-1\\\gcd(n,p-1)=1}} e^{\frac{i2\pi  \tau^n}{p}} \\
		&=&     \frac{1}{p} \sum_{1 \leq t\leq p-1} \left ( \sum_{1 \leq s\leq p-1} \omega^{-ts}e^{\frac{i2\pi  b \tau^s}{p}}-\sum_{1 \leq s\leq p-1} \omega^{-ts}e^{\frac{i2\pi  \tau^s}{p}}\right ) \nonumber \\
		&& \times \left (\sum_{d \mid p-1} \mu(d) \sum_{\substack{n \leq p-1, \\ d \mid n}}   \omega^{tn} \right ) \nonumber.
	\end{eqnarray}
	Taking the finite Fourier transform of the difference $D(b)=V(b)-V(1)$ returns 

	\begin{eqnarray} \label{eq9933RPI.500j}
		\widehat{D(a)}&=&	\sum_{1\leq b\leq  p-1}	 e^{-i2\pi \frac{ab}{p}}\left( \sum_{\substack{1\leq n\leq  p-1\\\gcd(n,p-1)=1}} e^{\frac{i2\pi b \tau^n}{p}} -\sum_{\substack{1\leq n\leq  p-1\\\gcd(n,p-1)=1}} e^{\frac{i2\pi  \tau^n}{p}}\right)  \\
		&=&  \frac{1}{p} \sum_{1\leq b\leq  p-1}	 e^{-i2\pi \frac{b}{p}}  \sum_{1 \leq t\leq p-1} \left ( \sum_{1 \leq s\leq p-1} \omega^{-ts}e^{\frac{i2\pi b \tau^s}{p}}-\sum_{1 \leq s\leq p-1} \omega^{-ts}e^{\frac{i2\pi a \tau^s}{p}}\right ) \nonumber \\
		&&\hskip 1.75in \times \left (\sum_{d \mid p-1} \mu(d) \sum_{\substack{n \leq p-1, \\ d \mid n}}   \omega^{tn} \right ) \nonumber\\
		&=&   \frac{1}{p} \sum_{1 \leq t\leq p-1} \left ( \sum_{1 \leq s\leq p-1} \omega^{-ts}\sum_{1\leq b\leq  p-1}	  e^{\frac{i2\pi b (\tau^s-a)}{p}}\right.  \nonumber \\
		&&\hskip .15in-\left .\sum_{1\leq b\leq  p-1}	 e^{-i2\pi \frac{ab}{p}}   \sum_{1 \leq s\leq p-1} \omega^{-ts}e^{\frac{i2\pi  \tau^s}{p}}\right ) \times \left (\sum_{d \mid p-1} \mu(d) \sum_{\substack{n \leq p-1, \\ d \mid n}}   \omega^{tn} \right ) \nonumber.
	\end{eqnarray}
	Now in the range $a<x=o(p)$, $\tau^s-a\ne0$ for any $s\in\{s<p-1:\gcd(s,p-1)=1\}$. Thus, using the geometric sum identity $\sum_{1\leq u\leq  p-1}	 e^{i2\pi au/p}=-1$ to simplify the last expression yields
	\begin{eqnarray} \label{eq9933RPI.500l}
		\widehat{D(a)}&=&	\sum_{1\leq b\leq  p-1}	 e^{-i2\pi \frac{ab}{p}}\left( \sum_{\substack{1\leq n\leq  p-1\\\gcd(n,p-1)=1}} e^{\frac{i2\pi b \tau^n}{p}} -\sum_{\substack{1\leq n\leq  p-1\\\gcd(n,p-1)=1}} e^{\frac{i2\pi  \tau^n}{p}}\right)  \nonumber\\
		&=&   \frac{1}{p} \sum_{1 \leq t\leq p-1} \left ( (-1)(-1)-(-1)  \sum_{1 \leq s\leq p-1} \omega^{-ts}e^{\frac{i2\pi  \tau^s}{p}}\right ) \nonumber \\
		&& \times \left (\sum_{d \mid p-1} \mu(d) \sum_{\substack{n \leq p-1, \\ d \mid n}}   \omega^{tn} \right ) .
	\end{eqnarray}
	Rearranging the last equation yield
	\begin{eqnarray} \label{eq9933RPI.500k}
		\widehat{V(a)}&=&\sum_{1\leq b\leq  p-1}	 e^{-i2\pi \frac{ab}{p}}\sum_{\substack{1\leq n\leq  p-1\\\gcd(n,p-1)=1}} e^{\frac{i2\pi b \tau^n}{p}} \\
		&=& -\sum_{\substack{1\leq n\leq  p-1\\\gcd(n,p-1)=1}} e^{\frac{i2\pi  \tau^n}{p}}  +  \frac{1}{p} \sum_{1 \leq t\leq p-1} \left ( 1-  \sum_{1 \leq s\leq p-1} \omega^{-ts}e^{\frac{i2\pi  \tau^s}{p}}\right ) \nonumber \\
		&&\hskip 2.5in \times \left (\sum_{d \mid p-1} \mu(d) \sum_{\substack{n \leq p-1, \\ d \mid n}}   \omega^{tn} \right ) \nonumber.
	\end{eqnarray}
	
	By \hyperlink{lem5555.400B}{Lemma} \ref{lem5555.400B}, the relatively prime summation kernel is bounded by
	\begin{eqnarray} \label{eq9933RPI.500m}
		\sum_{1 \leq t\leq p-1}	\Bigg |\sum_{d \mid p-1} \mu(d) \sum_{\substack{n \leq p-1, \\ d \mid n}}   \omega^{tn} \Bigg | 
		&=& \sum_{1 \leq t\leq p-1}\Bigg | \sum_{\gcd(n, p-1)=1}\omega^{tn} \Bigg |  \\ 
		&\ll &  p^{1+\delta}\log p\nonumber, 
	\end{eqnarray}
	where $\delta>0$ is a small number and by \hyperlink{lem1234A.150A}{Lemma} \ref{lem1234A.150A}, the difference including Gauss sum is bounded by
	\begin{eqnarray} \label{eq9933RPI.500o}
		\Bigg | 1-  \sum_{1 \leq s\leq p-1} \omega^{-ts}e^{\frac{i2\pi a \tau^s}{p}}\Bigg |=	\Bigg | 1- \sum_{1 \leq s\leq p-1} \chi(s) \psi(s) \Bigg| 
		&\leq & 2 p^{1/2} \log p, 
	\end{eqnarray}
	where  $\chi(s)=e^{i \pi s t/p}$, and $ \psi(s)=e^{i2\pi a \tau^s/p}$. Taking absolute value of the remainder term
	in (\ref{eq9933RPI.500k}) and replacing (\ref{eq9933RPI.500m}), and  (\ref{eq9933RPI.500o}), return
	\begin{eqnarray} \label{eq9933RPI.500p}
		|\widehat{R(a)}|	&=&\frac{1}{p} \Bigg |\sum_{1 \leq t\leq p-1} \Bigg ( 1-  \sum_{1 \leq s\leq p-1} \omega^{-ts}e^{\frac{i2\pi  \tau^s}{p}}\Bigg ) \cdot\Bigg (\sum_{d \mid p-1} \mu(d) \sum_{\substack{n \leq p-1, \\ d \mid n}}   \omega^{tn} \Bigg ) \Bigg | \nonumber\\
		&=&\frac{1}{p}\sum_{1 \leq t\leq p-1}\Bigg | 1-  \sum_{1 \leq s\leq p-1} \omega^{-ts}e^{\frac{i2\pi  \tau^s}{p}} \Bigg | \cdot \Bigg | \sum_{d \mid p-1} \mu(d) \sum_{\substack{n \leq p-1, \\ d \mid n}}   \omega^{tn} \Bigg | \nonumber\\
		&\ll &\frac{1}{p}(2 p^{1/2} \log p)\cdot \sum_{1 \leq t\leq p-1}\Bigg | \sum_{d \mid p-1} \mu(d) \sum_{\substack{n \leq p-1, \\ d \mid n}}   \omega^{tn} \Bigg | \nonumber\\
		&\ll &\frac{1}{p}(2 p^{1/2} \log p)\cdot (p^{1+\delta} \log p) \nonumber\\[.3cm]
		&\ll & p^{1/2+\delta} (\log p)^2,
	\end{eqnarray}
	where the implied constant depends on the number of divisors of $p-1$.
\end{proof}


\begin{lem}   \label{lem9933ERP.230V}\hypertarget{lem9933ERP.230V}  Let \(p\geq 2\) be a large prime. If $\tau $ be a primitive root modulo $p$ and $a<x=o(p)$ is not a primitive root, then
	\begin{equation} 
	\Bigg|\widehat{V(a)}\Bigg|=	\Bigg|\sum_{1\leq b\leq  p-1}	 e^{-i2\pi \frac{ab}{p}}	\sum_{\substack{1\leq n\leq p-1\\\gcd(n,p-1)=1}} e^{\frac{i2\pi b \tau^n}{p}}\Bigg| = O(p^{1/2+\delta} (\log p)^2)\nonumber,
	\end{equation} 
where $\delta>0$ is a small number and the implied constant is independent of $ b \in [1, p-1]$. 	
\end{lem} 
\begin{proof}[\textbf{Proof}] The second line in the estimation of the upper bound in \eqref{eq9933ERP.230d} follows from \hyperlink{thm9933ERP.220V}{Theorem} \ref{thm9933ERP.220V} and the fourth line follows from \hyperlink{thm9933Q.346}{Theorem} \ref{thm9933Q.346}: 
\begin{eqnarray} \label{eq9933ERP.230d}
\bigg|\widehat {V(u)}\bigg |&=& 	\Bigg|\sum_{1\leq b\leq  p-1}	 e^{-i2\pi \frac{ab}{p}}	\sum_{\substack{1\leq n\leq p-1\\\gcd(n,p-1)=1}} e^{\frac{i2\pi b \tau^n}{p}}\Bigg| \\[.2cm]
&=&\left |- \sum_{\substack{1\leq n\leq p-1\\\gcd(n,p-1)=1}} e^{i2\pi \frac{\tau ^n}{p}}+ O(p^{1/2+\delta} \log^2 p)  \right | \nonumber\\[.2cm]
&\ll &\left |\sum_{\substack{1\leq n\leq p-1\\\gcd(n,p-1)=1}} e^{i2 \pi \frac{\tau ^n}{p}} \right |+p^{1/2+\delta} (\log p)^2  \nonumber\\[.4cm]
&\ll&  p^{1/2+\delta} (\log p)^2\nonumber,
\end{eqnarray}
where $\delta>0$ is a small number. 
\end{proof}

\section{Evaluation Of The Main Term} \label{S9900-MT}\hypertarget{S9900-MT}
Finite sums and products over the primes numbers occur on various problems concerned with primitive roots. These sums and products often involve the normalized totient function $\varphi(n)/n=\prod_{p\mid n}(1-1/p)$ and the corresponding estimates, and the asymptotic formulas.

\begin{lemma} {\normalfont ({\color{red}\cite[Lemma 5]{MP1999}})}  \label{lem4.1}
		  Let \(x\geq 1\) be a large number, and let \(\varphi (n)\) be the Euler totient function.
		If \(q\leq \log^c x\), with $c\geq 0$ constant, an integer $1\leq a< q$ such that $\gcd(a,q)=1$, then
		\begin{equation}
		\sum_{\substack{p\leq x  \\ p \equiv a \bmod q} }\frac{\varphi(p-1)}{p-1}
		=A_q\frac{\li(x)}{\varphi(q)}
		+
		O\left(\frac{x}{\log ^bx}\right) ,\nonumber
		\end{equation}
		where \(\li(x)\) is the logarithm integral, and $b=b(c)>1$ is a constant, as \(x \rightarrow \infty\), and 
\begin{equation} \label{555}
A_q=\prod_{p \vert\gcd(a-1,q) } \left(1-\frac{1}{p}\right)\prod_{p \nmid q } \left(1-\frac{1}{p(p-1)}\right).
\end{equation} 
\end{lemma}

Related discussions for $q=2$ are given in {\color{red}\cite[Lemma 1]{SP1969}}, {\color{red}\cite[p.\ 16]{MP2004}}, and, \cite{VR1973}. The case \(q=2\) is ubiquitous in various results in Number Theory. 

\begin{lemma} \label{lem4.2}\hypertarget{lem4.2} Let \(x\geq 1\) be a large number, and let \(\varphi (n)\) be the Euler totient function.
		If \(q\leq \log^c x\), with $c\geq 0$ constant, an integer $1\leq a< q$ such that $\gcd(a,q)=1$, then
		\begin{equation} \label{el8859}
		\sum_{\substack{p\leq x \\ p \equiv a \bmod q}} \frac{1}{p}\sum_{\gcd(n,p-1)=1} 1=A_q \frac{\li(x)}{\varphi(q)}+O\left(\frac{x}{\log
			^b x}\right) ,\nonumber
		\end{equation} 
		where \(\li(x)\) is the logarithm integral, and $b=b(c)>1$ is a constant, as \(x \rightarrow \infty\), and $A_q$ is defined in {\normalfont (\ref{555})}. 
\end{lemma}

\begin{proof}[\textbf{Proof}] A routine rearrangement gives 
	\begin{eqnarray} \label{el500}
	\sum_{\substack{p\leq x \\ p \equiv a \bmod q}} \frac{1}{p}\sum_{\gcd(n,p-1)=1} 1&=&\sum_{\substack{p\leq x \\ p \equiv a \bmod q}} \frac{\varphi(p-1)}{p} \\[.3cm]
	&=&\sum_{\substack{p\leq x \\ p \equiv a \bmod q}} \frac{\varphi(p-1)}{p-1}-\sum_{\substack{p\leq x \\ p \equiv a \bmod q}} \frac{\varphi(p-1)}{p(p-1)} \nonumber.
	\end{eqnarray} 
	To proceed, apply Lemma \ref{lem4.1} to reach
	\begin{eqnarray} \label{el501}
	\sum_{\substack{p\leq x \\ p \equiv a \bmod q}}\frac{\varphi(p-1)}{p-1}  -\sum_{\substack{p\leq x \\ p \equiv a \bmod q}} \frac{\varphi(p-1)}{p(p-1)} 
	&=&A_q \frac{\li(x)}{\varphi(q)}+O \left (\frac{x}{\log^b x}\right )  \\
& & -\sum_{\substack{p\leq x \\ p \equiv a \bmod q}} \frac{\varphi(p-1)}{p(p-1)} \nonumber\\[.3cm]
	&=&A_q \frac{\li(x)}{\varphi(q)}+O \left (\frac{x}{\log^b x}\right ) \nonumber ,
	\end{eqnarray} 
	where the second finite sum 
\begin{equation}
\sum_{\substack{p\leq x \\ p \equiv a \bmod q}} \frac{\varphi(p-1)}{p(p-1)} \ll \log \log x
\end{equation} 
is absorbed into the error term,  $b=b(c)>1$ is a constant, and $A_q$ is defined in {\normalfont (\ref{555})}.   
\end{proof}


\section{Estimate For The Error Term} \label{S9900-ET}\hypertarget{S9900-ET}
The upper bound for exponential sum over subsets of elements in finite fields $\mathbb{F}_p$ stated in the last section will be used here to estimate the error term $E(x)$ arising in the proof of Theorem \hyperlink{thm9900.100}{Theorem} \ref{thm9900.100}. 

\begin{lemma}  \label{lem5.1}\hypertarget{lem5.1} Let \(p\geq 2\) be a large prime, let \(\psi \neq 1\) be an additive character, and let \(\tau\) be a primitive root mod \(p\). If the element \(u\ne 0\) is not a primitive root, then, 
\begin{equation} \label{el88400}
 \left \vert\sum_{\substack{x \leq p\leq 2x \\ p \equiv a \bmod q}}
\frac{1}{p}\sum_{\gcd(n,p-1)=1,} \sum_{ 0<k\leq p-1} \psi \left((\tau ^n-u)k\right) \right \vert\ll  x^{1/2+\delta}\log x,\nonumber 
\end{equation} 
where $1 \leq a <q, \, \gcd(a,q)=1$ and $O(\log^c x)$ with $c>0$ constant, for all sufficiently large numbers $x\geq 1$ and an arbitrarily small number \(\varepsilon<1/16\).
\end{lemma}

\begin{proof}[\textbf{Proof}] Let $\psi(z)=e^{i 2 \pi kz/p}$ with $0< k<p$, and rearrange the triple finite sum in the form
\begin{eqnarray} \label{e88901}
E(x)&=&  \sum_{\substack{x \leq p\leq 2x \\ p \equiv a \bmod q}}\frac{1}{p} \sum_{ 0<k\leq p-1,}  \sum_{\gcd(n,p-1)=1} \psi ((\tau ^n-u)k)  \\ 
&=&   \sum_{\substack{x \leq p\leq 2x \\ p \equiv a \bmod q}}
 \frac{1}{p}\sum_{ 0<k\leq p-1} e^{-i 2 \pi uk/p}   \sum_{\gcd(n,p-1)=1} e^{i 2 \pi k\tau ^n/p}\nonumber.
\end{eqnarray} 
Applying \hyperlink{lem9933ERP.230V}{Lemma} \ref{lem9933ERP.230V} to the inner double sum 
and completing the upper bound of the finite sum return
\begin{eqnarray} \label{el89991}
|E(x)|
	&\ll &  \sum_{\substack{x \leq p\leq 2x \\ p \equiv a \bmod q}} \frac{ 1}{p}  \cdot p^{1/2+\delta} (\log p)^2 \\[.2cm]
	&\ll & \sum_{x \leq p \leq 2x} \frac{ 1}{p}  \cdot p^{1/2+\delta} (\log p)^2 \nonumber\\[.2cm]
	&\ll &  \frac{ (\log x)^2}{x^{1/2-\delta }}\sum_{x \leq p \leq 2x}  1 \nonumber \\[.3cm]
	&\ll &  x^{1/2+\delta}\log x \nonumber,
\end{eqnarray}
where the number of primes in the short interval $[x,2x]$ is $\pi(2x)-\pi(x)\leq 2 x/ \log x$.

\end{proof}

A sharper estimate of the last finite sum over the primes in arithmetic progression can be estimated using the Brun-Titchmarsh theorem; this result states that the number of primes $p=qn+a$ in the interval $[x,2x]$ satisfies the inequality
\begin{equation} \label{el8040}
\pi(2x,q,a)-\pi(x,q,a) \leq \frac{3}{\varphi(q)}\frac{x}{ \log x},
\end{equation}
see {\color{red}\cite[p.\  167]{IK2004}}, {\color{red}\cite[p.\  157]{HG2007}}, \cite{MJ2012}, and {\color{red}\cite[p.\  83]{TG2015}}. \\

\section{The Main Result} \label{S9900-T1}\hypertarget{S9900-T1}
Given a fixed squarefree integer $u \ne \pm 1$, the precise primes counting function is defined by 
\begin{equation} \label{el8880}
\pi_{u}(x,q,a)=\# \{ p \leq x:p \equiv a \bmod q \text{ and } \ord_p(u)=p-1 \}
\end{equation}
for $1 \leq a <q$ and $\gcd(a,q)=1$. The limit
\begin{equation} \label{el8888}
\delta(u,q,a)=\lim_{x \to \infty} \frac{\pi_{u}(x,q,a)}{\pi(x,q,a)}=a_{u} \frac{A_q}{\varphi(q)}
\end{equation}
is the density of the subset of primes with a fixed squarefree primitive root $u\ne \pm 1$. 
\begin{theorem} \label{thm6.80} {\normalfont(\cite{HC1967}, \cite{LH1977})} Suppose the GRH is true. Then,
\begin{equation} \label{el6.80}
\pi_{u}(x,q,a)=\delta(u,q,a) \frac{x}{ \log  x}+O\left(\frac{x \log \log x}{ \log^2  x}\right).\nonumber
\end{equation}
\end{theorem}
As explained in {\color{red}\cite[Section 8.1]{MP2004}}, 
the existing primitive roots counting method fails to prove any unconditional result on primes and primitive roots. To circumvent this 
obstacle, the proof of \hyperlink{thm9900.100}{Theorem} \ref{thm9900.100} below, uses a new primitive roots counting method explicated in \hyperlink{S9955D}{Section} \ref{S9955D}. \\

\begin{proof}[\textbf{Proof}] (\hyperlink{thm9900.100}{Theorem} \ref{thm9900.100}) Suppose that the squarefree integer $u \ne \pm 1$ is not a primitive root for all primes \(p\geq x_0\), with \(x_0\geq 1\) constant. Let \(x>x_0\) be a large number, and $q=O(\log^cx)$. Consider the sum of the characteristic function over the short interval \([x,2x]\), that is, 
\begin{equation} \label{el8720}
	\pi_{u}(x,q,a)=\sum_{\substack{x \leq p\leq 2x \\ q \equiv a \bmod q}} \Psi (u)=0.
\end{equation}
Replacing the characteristic function, \hyperlink{lem9955.200A}{Lemma} \ref{lem9955.200A}, and expanding the nonexistence equation (\ref{el8720}) yield
\begin{eqnarray} \label{el8730}
\pi_{u}(x,q,a)&=&\sum _{\substack{x \leq p\leq 2x\\
			p \equiv a \bmod q		}} \Psi (u)  \\ [.2cm]
	&=&\sum_{\substack{x \leq p\leq 2x\\[.2cm]
			p \equiv a \bmod q		}} \left (\frac{1}{p}\sum_{\gcd(n,p-1)=1,} \sum_{ 0\leq k\leq p-1} \psi \left((\tau ^n-u)k\right) \right ) \nonumber\\[.2cm]
	&=&a_{u}\sum_{\substack{x \leq p\leq 2x\\
			p \equiv a \bmod q		}} \frac{1}{p}\sum_{\gcd(n,p-1)=1} 1 \nonumber\\[.2cm]
&& \hskip 1.5 in +\sum_{\substack{x \leq p\leq 2x\\
			p \equiv a \bmod q		}}
	\frac{1}{p}\sum_{\gcd(n,p-1)=1,} \sum_{ 0<k\leq p-1} \psi \left((\tau ^n-u)k\right)\nonumber\\[.3cm]
	&=&M(x) + E(x)\nonumber,
	\end{eqnarray} 
	where $a_{u} \geq 0$ is a constant depending on the integers $u\ne \pm 1$ and $q\geq 2$. \\
	
The main term $M(x)$ is determined by a finite sum over the trivial additive character \(\psi =1\), and the error term $E(x)$ is determined by a finite sum over the nontrivial additive characters \(\psi(z) =e^{i 2\pi  z/p}\neq 1\).\\
	
	Take a constant $b=b(c)>1$, depending on $c\geq0$. Applying \hyperlink{lem4.2}{Lemma} \ref{lem4.2} to the main term, and \hyperlink{lem5.1}{Lemma} \ref{lem5.1} to the error term yield
	\begin{eqnarray} \label{el8762}
\sum _{\substack{x \leq p\leq 2x\\
			p \equiv a \bmod q		}} \Psi (u) 
	&=&M(x) + E(x) \\[.2cm]
	&=&a_{u}\left (A_q\frac{\li(2x)-\li(x)}{\varphi(q)} \right )+O\left(\frac{x}{\log^bx}\right)+O\left(x^{1/2+\delta}\log x \right) \nonumber\\[.4cm]
	&=& \delta(u,q,a)\left (\li(2x)-\li(x)\right )+O\left( \frac{x }{\log^b} \right)  \nonumber,
	\end{eqnarray} 
	where $\delta(u,q,a)=a_{u}A_q /\varphi(q) \geq 0$, and $a_{u}\geq 0$ is a correction factor depending on the squarefree integer $u\ne\pm1$. \\

But $\delta(u,q,a) > 0$ contradicts the hypothesis \eqref{el8720} for all sufficiently large numbers $x \geq x_0$. Ergo, the short interval $[x,2x]$ contains primes $p=qn+a$ such that the $u$ is a fixed primitive root. Specifically, the counting function is given by
	\begin{equation} \label{el8964}
	\pi_u(x,q,a)=\sum _{\substack{p\leq x\\
			p \equiv a \bmod q		}} \Psi (u)
	=\delta(u,q,a)\li(x) +O\left( \frac{x}{\log^b x} \right) .
	\end{equation} 
	This completes the verification. 
\end{proof}

The determination of the correction factor $a_{u}$ in a primes counting problem is a complex problem, some cases are discussed in \cite{SP2003}, and \cite{ZD2009}. \\


\end{document}